\documentclass{hha}

\usepackage[all]{xy}
\usepackage{enumerate}

\newtheorem{thmA}{Theorem}

\theoremstyle{plain}
\newtheorem{conjA}{Conjecture}

\newcommand{\F}{\mathbb{F}}
\newcommand{\Zp}{\mathbb{Z}_{(p)}}
\newcommand{\Z}{\mathbb{Z}}
\newcommand{\Mod}{\text{-Mod}}

\newcommand{\Fc}{\mathcal{F}}
\newcommand{\Li}{\mathcal{L}}
\newcommand{\T}{\mathcal{T}}
\newcommand{\Bc}{\mathcal{B}}
\newcommand{\Cc}{\mathcal{C}}
\newcommand{\Ec}{\mathcal{E}}
\newcommand{\h}{\mathcal{H}}

\newcommand{\pcompl}[1]{#1^{\wedge}_p}
\newcommand{\limproj}[1]{\lim\limits_{\substack{\longleftarrow \\ #1}}}

\newcommand{\Ob}{\text{Ob}}
\newcommand{\Hom}{\text{Hom}}
\newcommand{\Aut}{\text{Aut}}
\newcommand{\Out}{\text{Out}}
\newcommand{\Inn}{\text{Inn}}
\newcommand{\Inj}{\text{Inj}}
\newcommand{\Mor}{\text{Mor}}
\newcommand{\Res}{\text{Res}}
\newcommand{\Id}{\text{Id}}
\newcommand{\im}{\text{Im}\,}
\newcommand{\Ind}{\text{Ind}}
\newcommand{\coInd}{\text{coInd}}
\newcommand{\Syl}{\text{Syl}}

\newcommand{\hyp}{\mathfrak{hyp}}

\begin{document}

\title{Cohomology of linking systems with twisted coefficients by a $p$-solvable action} 

\shorttitle{Cohomology of linking systems with $p$-solvable action}

\author{R\'emi	Molinier}             

\email{molinier@ksu.edu}

%
%
\address{Department of Mathematics,
         Kansas State University,
         138 Cardwell Hall,
         Manhattan, KS 66506,
         USA}







\begin{abstract}
In this paper we study the cohomology of the geometric realization of linking systems with twisted coefficients. More precisely, given a prime $p$ and a $p$-local finite group $(S,\Fc,\Li)$, we compare the cohomology of $\Li$ with twisted coefficients with the submodule of $\Fc^c$-stable elements in the cohomology of $S$. We start with the particular case of constrained fusion systems. Then, we study the case of $p$-solvable actions on the coefficients. 
\end{abstract}

\received{Month Day, Year}   
\revised{Month Day, Year}    
\published{Month Day, Year}  
\submitted{Nathalie Wahl}      
\volumeyear{} 
\volumenumber{} 
\issuenumber{}   
\startpage{1}     
\articlenumber{} 
\owner{International Press}

\maketitle


\section{Introduction}

The notion of saturated fusion system was introduced by Puig in the 90s in a context of modular representation theory. Since their introduction, topologists use them to study classifying spaces of finite groups or, more precisely, their $p$-completions.
A \emph{$p$-local finite group} is a triple $(S,\Fc,\Li)$ where $S$ is a $p$-group, $\Fc$ a saturated fusion system over $S$ and $\Li$ an associated linking system.
For a $p$-local finite group $(S,\Fc,\Li)$, $\pcompl{|\Li|}$ is called its \emph{classifying space}. 
The theory of $p$-local finite groups have been studied in details by Broto, Levi, Oliver and others 
(see \cite{BLO2}, \cite{OV1}, \cite{5a1} and \cite{5a2}). 
The linking system and its geometric realization, even without $p$-completion, play here a fundamental and central role. 
In fact, for a given saturated fusion systems, the existence and uniqueness of an associated linking system were shown more recently by Chermak \cite{Ch} (using the theory of partial groups). The proof of this important conjecture highlights that linking systems and their geometric realizations form a deep link between fusion system theory and homotopy theory (we refer to Aschbacher, Kessar and Oliver \cite{AKO} for more details about fusion systems in general).

An old and well-known result due to Cartan and Eilenberg (see \cite[Theorem XII.10.1]{CE}) 
expresses the cohomology of a finite group $G$ in a $\Zp[G]$-module as the submodule 
of ``stable'' elements in the cohomology of a Sylow $p$-subgroup of $G$.
This submodule of stable elements corresponds to the inverse limit over the ``fusion`` of the group cohomology functor.
One important result in the theory of $p$-local finite groups is an analog of this theorem for $p$-local finite groups which tells us that the cohomology of the geometric realization of a linking system can be computed by $\Fc$-stable elements. More precisely, there is a natural inclusion of $BS$ into $|\Li|$ and it induces the following isomorphism. Here, $\Fc^c$ is the full subcategory of $\Fc$ consisting of $\Fc$-centric subgroups of $S$ and, for $A$ a finite $\Zp$-module, $H^*(\Fc^c,A)\subseteq H^*(S,A)$ is the submodule of $\Fc$-stable elements.

\begin{theorem}\label{trivial}
Let $(S,\Fc,\Li)$ be a $p$-local finite group and $A$ be a finite $\Zp$-module.
The inclusion of $BS$ in $|\Li|$ induces a natural isomorphism
\[
 \xymatrix{ H^*(\pcompl{|\Li|},A)\cong H^*(|\Li|,A)\ar[r]^-\cong & H^*(\Fc^c,A).}
\]
\end{theorem}

\begin{proof}
The case $A=\F_p$ is \cite[Theorem B]{BLO2} and the general case is proven in \cite[Lemma 6.12]{5a2}.
\end{proof}

One question asked by Oliver in his book with Aschbacher and  Kessar \cite{AKO} is the understanding of the cohomology of $|\Li|$ with twisted coefficients. 
Indeed, this cohomology appears for example in the study of extensions of $p$-local finite groups or, more directly, can give more information about the link between the fusion system and the spaces $|\Li|$ or $\pcompl{|\Li|}$. 
Recall that, if a space $X$ has a universal covering space $\widetilde{X}$, the cohomology of $X$ with twisted coefficients in a $\Z[\pi_1(X)]$-module $M$ is the cohomology of the chain complex
\[C^*(X;M)=\Hom_{\Z[\pi_1(X)]}(S_*(\widetilde{X}),M),\]
where $S_*(\widetilde{X})$ is the usual singular chain complex of $\widetilde{X}$.

Levi and Ragnarsson \cite{LR} already give some tools along these lines. In an other paper \cite{Mo1}, the author extends Theorem \ref{trivial} to the case of nilpotent actions on the coefficients. The main ingredient  is to construct, as in the trivial coefficient case, an idempotent of $H^*(S,M)$ with image $H^*(\Fc^c,M)$.

In this paper, we also want to extend Theorem \ref{trivial} to twisted coefficients but when the action factors through a $p$-solvable group. The methods used here are completely different from the ones used in \cite{Mo1} and also more direct. We first have a look at constrained fusion systems. In that case we are able to prove that, with any coefficient module, the cohomology of $|\Li|$ can be computed by stable elements.

\begin{thmA}[see Theorem \ref{isoconstrained}]
Let $(S,\Fc,\Li)$ be a $p$-local finite group.
 If $\Fc$ is constrained and $M$ is a $\Zp[\pi_1(|\Li|)]$-module, then the inclusion of $BS$ in $|\Li|$ induces an isomorphism,
 \[
  H^*(|\Li|,M)\cong H^*(\Fc^c,M).
 \]
\end{thmA}

Next we focus on $p$-solvable actions. The main ingredients here are $p$-local finite subgroups of index a power of $p$ or prime to $p$ and their homotopy properties. We start by looking at $p$-local subgroups of index prime to $p$ (see Definition \ref{p and p' sub}(b)).

\begin{thmA}[see Theorem \ref{main}]
Let $(S,\Fc,\Li)$ be a $p$-local finite group and denote by $(S,O^{p'}(\Fc),O^{p'}(\Li))$  its minimal $p$-local subgroup of index prime to $p$.
If $M$ is a $\Zp[\pi_1(|\Li|)]$-module and if the inclusion of $BS$ in $|O^{p'}(\Li)|$ induces an isomorphism 
\[H^*(|O^{p'}(\Li)|,M)\cong H^*(O^{p'}(\Fc)^c,M),\] 
then the inclusion of $BS$ in $|\Li|$ induces an isomorphism
\[
 H^*(|\Li|,M)\cong H^*(\Fc^c,M).
\]
\end{thmA} 

This Theorem allows us to prove that if the action on the coefficients factor through a $p'$-group or, even better, a $p$-nilpotent group, then the cohomology of $|\Li|$ can be computed by stable elements. 

It is much more complicated to work with $p$-local finite groups of index a power of $p$, especially on the level of stable elements. Indeed, for $(S_0,\Fc_0,\Li_0)$ a $p$-local subgroup of $(S,\Fc,\Li)$ of index a power of $p$ and $M$ a $\Zp[\pi_1(|\Li|)]$-module, it is difficult to compare $H^*(\Fc^c,M)$ and $H^*(\Fc_0^c,M)$. The difficulty mostly comes from the fact that we are working on different $p$-groups: $S$ and $S_0$. But when we work with a $p$-local finite group realizable by a finite group $G$, and if $G$ acts "consistently" on the coefficients it is possible to get some positive results (see Section \ref{sec psolvable}).

\begin{thmA}[see Corollary \ref{psolvable2}] 
 Let $G$ be a finite group, $S$ a Sylow $p$-subgroup of $G$ and $(S,\Fc,\Li)$ the associated $p$-local finite group.
 Let $M$ be a $\Zp[\pi_1(|\Li|)]$-module and assume that $G$ acts \emph{consistently} on $M$.
If both actions factor through a given $p$-solvable $\Gamma$ and all the $M$-essential subgroups (see Definition \ref{M essential}) of $S$ are $p$-centric, then we have natural isomorphisms,
\[H^*(|\Li|,M)\cong H^*(G,M) \cong H^*(\Fc^c,M).\]
\end{thmA}

All of these results lead us to the following conjecture.

\begin{conjA}[see Conjecture \ref{conjecture}]
Let $(S,\Fc,\Li)$ be a $p$-local finite group and $M$ a $\Zp[\pi_1(|\Li|)]$-module.
 If the action of $\pi_1(|\Li|)$ on $M$ is $p$-solvable, then the inclusion of $BS$ in $|\Li|$ induces a natural isomorphism
\[
 \xymatrix{ H^*(|\Li|,M)\ar[r]^\cong & H^*(\Fc^c,M).}
\]  
\end{conjA}

We finish this paper with an example for Conjecture \ref{conjecture} which does not follow from the the other results.

\textbf{Organization.} In Section \ref{sec background} we give some background on $p$-local finite groups and stable elements. Section \ref{sec constrained} is dedicated to the case of constrained fusion systems, Section \ref{sec p'} to coprime actions  and Section \ref{sec psolvable} to $p$-solvable actions for a realizable $p$-local finite group. Finally, we give in Section \ref{sec example} an example for Conjecture \ref{conjecture}.

\ack
I would like to thanks Bob Oliver, my PhD adviser, for his help and support all along this work. 
I am also grateful to Jesper Grodal for his hospitality at the Center for Symmetry and Deformation and many fruitful conversations.
Finally, I would like to thanks the referee for his very careful reading of the paper and all his useful and accurate comments.

\section{Background}\label{sec background}

We give here a very short introduction to $p$-local finite groups. 
The notion of fusion system was first introduced by Puig for modular representation theory purpose.
Later, Broto, Levi and Oliver developed the notion of linking systems and $p$-local finite groups to study $p$-completed classifying spaces of finite groups and spaces which have similar homotopy properties. We refer the reader interested in more details to Aschbacher, Kessar and Oliver \cite{AKO}.

\subsection{Fusion systems and linking systems}
A fusion system over a $p$-group $S$ is a way to abstract the action of a finite group $G$ with $S\in\Syl_p(G)$ on the subgroups of $S$ by conjugation. 
For $G$ a finite group and $g\in G$, we will denote by $c_g$ the homomorphism $x\in G\mapsto gxg^{-1}\in G$ and for $H,K$ two subgroups of $G$, 
$\Hom_G(H,K)$ will denote the set of all group homomorphism $c_g$ for $g\in G$ such that $c_g(H)\leq K$.

\begin{definition}\label{defF}
 Let $S$ be a finite $p$-group.
 A \emph{fusion system} over $S$ is a small category $\Fc$, where $\Ob(\Fc)$ is the set of all subgroups of $S$ and which satisfies the following two
 properties for all $P,Q\leq S$:
 \begin{enumerate}[(a)]
  \item $\Hom_S(P,Q)\subseteq \Mor_\Fc(P,Q)\subseteq \Inj(P,Q)$;
  \item each $\varphi\in\Mor_\Fc(P,Q)$ is the composite of an $\Fc$-isomorphism followed by an inclusion.
 \end{enumerate}
A fusion system is \emph{saturated} if it satisfy two more technical axioms called the saturation axioms (we refer the reader to \cite[Definition I.2.1]{AKO} for a proper definition).
\end{definition}

The composition in a fusion system is given by composition of homomorphisms. We usually write $\Hom_\Fc(P,Q)=\Mor_\Fc(P,Q)$ to emphasize that the morphims in $\Fc$ are homomorphisms. For $P,Q\leq S$, we say that $P$ is \emph{$\Fc$-conjugate} to $Q$ if there is an $\Fc$-isomorphism between $P$ and $Q$. We denote by $P^\Fc$ the set of all subgroups of $S$ which are $\Fc$ conjugate to $P$.

The typical example of a saturated fusion system is the fusion system $\Fc_S(G)$ of a finite group $G$ over $S\in\Syl_p(G)$. 

For the purpose of this paper, we need to distinguish some collections of subgroups of $S$.

\begin{definition}
Let $\Fc$ be a saturated fusion system over a finite $p$-group $S$.
\begin{enumerate}[(a)]
\item A subgroup $P\leq S$ is \emph{$\Fc$-centric} if for every $Q\in P^\Fc$, $C_S(Q)=Z(Q)$.
\item A subgroup $P\leq S$ is \emph{$\Fc$-radical} if $O_p(\Aut_\Fc(P))=\Inn(P)$.
\item A subgroup $P\leq S$ is \emph{$\Fc$-quasicentric} if for each $Q\leq PC_S(P)$ containing $P$, and each $\alpha\in\Aut_\Fc(Q)$ such that $\alpha|_P=\Id$, $\alpha$ has a $p$-power order. 
\end{enumerate}
We let $\Fc^{cr}\subseteq \Fc^{c}\subseteq\Fc^{q}\subseteq \Fc$ denote the full subcategories of $\Fc$ with objects the $\Fc$-centric and $\Fc$-radical subgroups, the $\Fc$-centric subgroups  and the $\Fc$-quasicentric subgroups, respectively.
\end{definition}

If $\Fc=\Fc_S(G)$, a subgroup $P\leq S$ is 
\begin{enumerate}[(a)]
\item \emph{$\Fc$-centric} if and only if it is $p$-centric (i.e. $Z(P)\in\Syl_p(C_G(P))$, 
\item \emph{$\Fc$-radical} if $P/Z(P)=O_p(N_G(P)/C_G(P))$.
\item \emph{$\Fc$-quasicentric} if and only if $O^p(C_G(P))$ has order prime to $p$. 
\end{enumerate}

The notion of linking system has been introduced by Broto, Levi and Oliver \cite{BLO2} and generalized by Broto, Castellana, Grodal and Oliver in \cite{5a1}. We refer the reader to theses papers, or \cite[Part III]{AKO}, for a proper definition. We recall here some basic facts about linking systems which will be needed here.

For $G$ a finite group, $S\in\Syl_p(G)$ and $\h$ a collection of subgroups of $S$, the \emph{transporter category} of $G$ over $S$ with set of objects $\h$ is the category $\T_H^\h(G)$  with objects $\h$ and for $P,Q\in \h$, $\Mor_\Li(P,Q)=T_G(P,Q)=\{g\in G\mid P^g\leq Q\}$.
For $\Fc$ a saturated fusion system over a $p$-group $S$, a \emph{linking system} associated to $\Fc$ is a certain finite category with objects a collection $\h$ of subgroups of $S$ together with two functors:
\[\xymatrix{\T_S^\h(S)\ar[r]^-{\delta} & \Li \ar[r]^-\pi & \Fc.}\]
$\delta$ is the identity on objects and injective on morphisms and $\pi$ is injective on objects and surjective on morphisms. 
The collection $\h$ has to be stable by overgroups and $\Fc$-conjugation and the following proposition tell you which collection you can have.

\begin{proposition}\label{object linking system}
 Let $\Fc$ be a saturated fusion system over a $p$-group $S$.
 Let $\Li$ be a linking system associated to $\Fc$.
 \begin{enumerate}[(a)]
 \item $\Ob(\Fc^{cr})\subseteq \Ob(\Li)\subseteq \Ob(\Fc^q)$, and there exists a linking system $\Li^q$ associated to $\Fc$ such that $\Ob(\Li^q)=\Ob(\Fc^q)$, and $\Li$ is a full subcategory of $\Li^q$.
 
 \item For every subset $\Ob(\Fc^{cr})\subseteq\h\subseteq \Ob(\Fc^q)$ stable by $\Fc$-conjugacy and overgroups, the full subcategory $\Li^\h$ of $\Li^q$ with set of objects $\h$ is also a linking system associated to $\Fc$.
\end{enumerate}
\end{proposition}

\begin{proof}
 The first point of (a) can be found for example in \cite[Proposition 4(g)]{O4}. For the second statement of (a), you can find a proof in \cite[Proposition III.4.8]{AKO}.
 Finally, (b) is a consequence of the definition of linking systems.
\end{proof}

If $\h=\Ob(\Fc^q)$, $\Li$ is called a \emph{quasicentric linking system} and if $\h=\Ob(\Fc^c)$, $\Li$ is called a \emph{centric linking system}.

\begin{definition}
A $p$-local finite group is a triple $(S,\Fc,\Li)$ where $\Fc$ is a saturated fusion system over $S$ and $\Li$ is an associated linking system.
If $(S_0,\Fc_0,\Li_0)$ is an other $p$-local finite group, we will say that $(S_0,\Fc_0,\Li_0)$ is a \emph{$p$-local subgroup} of $(S,\Fc,\Li)$ if $S_0\leq S$ and $\Fc_0\subseteq \Fc$ is a subsystem of $\Fc$. Notice that we do not require that $\Li_0$ is a subcategory of $\Li$.
\end{definition}

The typical example you should have in mind is the following. For $G$ a finite group and $S\in\Syl_p(G)$ let $\Li_S^q(G)$ be the category with objects the $\Fc_S(G)$-quasicentric subgroups of $G$ and, for $P,Q\in\Ob(\Li)$,
 \[\Mor_\Li(P,Q)=T_G(P,Q)/O^{p}(C_G(P)).\] 
 Then $(S,\Fc_S(G),\Li_S^q(G))$  defines a $p$-local finite group where $\Li_S^q(G)$ is a quasicentric linking system. We also denote by $\Li_S^c(G)$ the full subcategory of $\Li_S^q(G)$ with objects the $p$-centric subgroups of $S$ and it is a centric linking system.

We finish with some basic homotopy properties about linking systems which will be needed in this paper. We refer the reader interested in more details to \cite[Part III]{AKO}. For $(S,\Fc,\Li)$ a $p$-local finite group, we write $|\Li|$ for the geometric realization of $\Li$ and $\pi_\Li=\pi_1(|\Li|)$ for its fundamental group. The following Theorem will allow us to change the set of objects of $\Li$ without changing the homotopy type of $|\Li|$.

\begin{theorem}[{\cite[Theorem 3.5]{5a1}}]\label{differentLisimeq}
Let $\Fc$ be a saturated fusion system over a $p$-group $S$.
Let $\Li_0\subseteq \Li$ be two linking systems associated to $\Fc$ with a different set of objects.
Then the inclusion induces a homotopy equivalence of space $|\Li_0|\simeq |\Li|$.
\end{theorem}

\subsection{$p$-local finite subgroups of index a power of $p$ or prime to $p$}

The notions $p$-local subgroups of index a power of $p$ or prime to $p$ have been introduced and studied by Broto, Castellana, Grodal, Levi and Oliver \cite{5a2}. 
Here we just give the definitions what we need about these $p$-local subgroups and we refer the reader to \cite{5a2} for more details.

\begin{definition}\label{p and p' sub}
 Let $(S,\Fc,\Li)$ be a $p$-local finite group and $(S_0,\Fc_0,\Li_0)$ a $p$-local subgroup of $(S,\Fc,\Li)$.
  Set $\hyp(\Fc)=\langle g^{-1}\alpha(g)\;\mid\; g\in P\leq S, \alpha\in O^p\left(\Aut_\Fc(P)\right)\rangle\trianglelefteq S$.
 \begin{enumerate}[(a)]
  \item We say that $(S_0,\Fc_0,\Li_0)$ is a \emph{$p$-local subgroup of index a power of $p$} if
       $S_0\geq \hyp(\Fc)$ and, for every $P\leq S_0$, $O^p(\Aut_\Fc(P))\leq \Aut_{\Fc_0}(P)$.
   
  \item We say that $(S_0,\Fc_0,\Li_0)$ is a \emph{$p$-local subgroup of index prime to $p$} if
       $S_0=S$ and, for every $P\leq S$, $O^{p'}(\Aut_\Fc(P))\leq \Aut_{\Fc_0}(P)$.
 \end{enumerate}
\end{definition}

Notice that $\hyp(\Fc)$ is denoted $O^p_\Fc(S)$ in \cite[Definition 2.1]{5a2}.
These particular $p$-local subgroups satisfy the following properties.

\begin{proposition}[{\cite[Proposition 3.8]{5a2}}]\label{p and p' sub prop}
 Let $(S,\Fc,\Li)$ be a $p$-local finite group and $(S_0,\Fc_0,\Li_0)$ a $p$-local subgroup of $(S,\Fc,\Li)$.
 \begin{enumerate}[(a)]
  \item If $(S_0,\Fc_0,\Li_0)$ is of index a power of $p$, then $P\leq S_0$ is $\Fc_0$-quasicentric if, and only if, $P$ is $\Fc$-quasicentric.
  \item If $(S_0,\Fc_0,\Li_0)$ is of index prime to $p$, then $P\leq S$ is $\Fc_0$-centric if, and only if, $P$ is $\Fc$-centric.
 \end{enumerate}
\end{proposition} 

For an infinite group $G$, we denote by $O^{p'}(G)$ the intersection of all normal subgroups in $G$ of finite index prime to $p$. For $\Fc$ a fusion system over a $p$-group $S$, let $O^{p'}_*(\Fc)$ be the fusion system generated by $O^{p'}(\Aut_\Fc(P))$ for all $P\leq S$
 and define
\[
  \Out_\Fc^0(S)=\langle \alpha\in\Out_\Fc(S)\;\mid\; \alpha|_P\in\Hom_{O^{p'}_*(\Fc)}(P,S)\text{, for some }P\leq S\rangle. 
\]
Since $\Aut_\Fc(S)$ normalizes $O^{p'}_*(\Fc)$, $\Out_\Fc^0(S)\trianglelefteq \Out_\Fc(S)$.

\begin{proposition}\label{p'generate}
 Let $(S,\Fc,\Li)$ be a $p$-local finite group.
 \begin{enumerate}[(a)]
  \item $\Fc=\langle \Aut_\Fc(S),O^{p'}_*(\Fc)\rangle$.
  \item $\pi$ and the inclusion of $B\Aut_\Fc(S)$ in $|\Fc^c|$ induce isomorphisms,
           \[
            \xymatrix{\theta\colon\pi_\Li/O^{p'}(\pi_\Li)\ar[r]^-{\cong} & \pi_1(|\Fc^c|)\ar[r]^-{\cong} & \Out_\Fc(S)/\Out_\Fc^0(S).}
           \]
 \end{enumerate}
\end{proposition}

\begin{proof}
 The point $(a)$ is proved in \cite[Lemma 3.4]{5a2}. For $(b)$, the second isomorphism is given in \cite[Proposition 5.2]{5a2} and the first one in \cite[Theorem 5.5]{5a2} and the comment which follows. 
\end{proof}

According to Proposition \ref{p and p' sub prop}, when dealing with $p$-local subgroups of index prime to $p$, we will work with centric linking systems.

\begin{theorem}[{\cite[Theorem 5.5]{5a2}}]\label{p'cover}
 Let $(S,\Fc,\Li)$ be a $p$-local finite group with $\Li$ a centric linking system.
 For each subgroup $H\leq\Out_\Fc(S)$ containing $\Out_\Fc^0(S)$, there is a unique $p$-local finite subgroup $(S,\Fc_H,\Li_H)$ of index prime to $p$ such that $\Out_{\Fc_H}(S)=H$ and $\Li_H=\pi^{-1}(\Fc_H^c)$.
 
 Moreover, $|\Li_H|$ is homotopy equivalent, via its inclusion in $|\Li|$, to the covering space of $|\Li|$ with fundamental group $\widetilde{H}\geq O^{p'}(\pi_\Li)$ such that $\theta(\widetilde{H}/O^{p'}(\pi_\Li))=H/\Out_\Fc^0(S)$ 
 (where  $\theta$ is the isomorphism given in Proposition \ref{p'generate}(b)).
\end{theorem}

Thus, for a $p$-local finite group $(S,\Fc,\Li)$, with $\Li$ a centric linking system, 
we can define the \textit{minimal $p$-local subgroup of index prime to $p$}, $(S,O^{p'}(\Fc),O^{p'}(\Li))$
corresponding to $(S,\Fc_H,\Li_H)$ with $H=\Out_\Fc^0(S)$ in Theorem \ref{p'cover}.

\subsection{Cohomology and stable elements}

The first result about stable elements is due to Cartan and Eilenberg (\cite[Chap XII, Theorem 10.1]{CE}). It also served as a guideline in the establishment of Theorem \ref{trivial} by Broto, Levi and Oliver. Here we recall the definition of $\Fc^c$-stable elements in a context of twisted coefficients. We refer the reader to \cite{Mo1} for more details. 
As in \cite{Mo1}, we will denote by $\omega\colon\Li\rightarrow \pi_\Li=\pi_1(|\Li|,S)$ the functor which maps each object to the unique object in the target and sends each morphism $\varphi\in\Mor_\Li(P,Q)$ to the class of the loop $\iota_Q.\varphi. \overline{\iota_P}$ where $\iota_P=\delta(1)\in\Mor_\Li(P,S)$, $\iota_Q=\delta(1)\in\Mor_\Li(Q,S)$ and  $\overline{\iota_P}$ is the edge $\iota_P$ followed in the opposite direction.

Let $(S,\Fc,\Li)$ be a $p$-local finite group. Recall first that $\delta\colon\T_S^{\Ob(\Li)}(S)\rightarrow \Li$ induces an inclusion $\delta_S\colon BS\rightarrow |\Li|$. In particular, it induces a natural map $S\rightarrow \pi_\Li$ and thus, for every $\Zp[\pi_\Li]$-module $M$, we have a natural action of $S$, or any subgroup of $S$, on $M$. 
Now, let $M$ be a $\Zp[\pi_\Li]$-module, the group cohomology bifunctor $H^*(-,-)$ induces a functor
\[\xymatrix{ H^*(-,M) \colon\Fc^c\ar[r] & \Zp\Mod}\]
(a priori, $H^*(g,M)$ is defined for $g\in\Mor(\Li)$ but \cite[Proposition 2.2]{Mo1} proves that $H^*(-,M)$ is well defined on $\Fc^c$).

\begin{definition}\label{fcstable}
Let $(S,\Fc,\Li)$ be a $p$-local finite group.
An element $x\in H^*(S,M)$ is called \emph{$\Fc$-centric stable}, or \emph{$\Fc^c$-stable}, if for all $P\in\Ob(\Fc^c)$ and all $\varphi\in\Hom_\Fc(P,S)$,
 \[\varphi^*(x)=\Res_P^S(x).\]
We denote by $H^*(\Fc^c,M)\subseteq H^*(S,M)$ the submodule of all $\Fc^c$-stable elements.
\end{definition}

Notice that
\[H^*(\Fc^c,M)=\limproj{\Fc^c} H^*(-,M)=\limproj{\Li} H^*(-,M)\]
where the last equality holds if $\Li$ is a centric linking system.

\section{Constrained fusion systems}\label{sec constrained}
Let $(S,\Fc,\Li)$ be a $p$-local finite group. 
Here, we assume that $\Fc$ is a constrained fusion system.

\begin{definition}\label{constrained}
Let $\Fc$ be a fusion system over a $p$-group $S$.
A subgroup $Q\leq S$ is \emph{normal in $\Fc$} if 
\begin{enumerate}[(i)]
\item $Q\trianglelefteq S$, and
\item for all $P,R\leq S$ and every $\varphi\in\Hom_\Fc(P,R)$, $\varphi$ extends to a morphism
$\overline{\varphi}\in\Hom_\Fc(PQ,RQ)$ such that $\overline{\varphi}(Q)=Q$. 
\end{enumerate} 
We write $O_p(\Fc)$ for the maximal subgroup of $S$ which is normal in $\Fc$.
We say that  $\Fc$ is \emph{constrained} if $O_p(\Fc)$ is $\Fc$-centric. 
\end{definition}

An important and classical result about constrained fusion systems is the following.
\begin{proposition}[{\cite[Proposition 4.3]{5a1}}]\label{model}
 Let $(S,\Fc,\Li)$ be a $p$-local finite group with $\Li$ a centric linking system.
 If $\Fc$ is constrained, there exists a finite group $G$ such that
 \begin{enumerate}[(a)]
  \item $S$ is a Sylow $p$-subgroup of $G$,
  \item $C_G(O_p(G))\leq O_p(G)$,
  \item $\Fc_S(G)=\Fc$.
 \end{enumerate}
Moreover, $G\cong\Aut_\Li(O_p(\Fc))$ and $\Li\cong\Li_S^c(G)$.
\end{proposition}
This group $G$ is called a \emph{model} of $\Fc$ and it is unique in a precise way (see \cite[Theorem III.5.10]{AKO}). This model can also be recovered from the homotopy type of the geometric realization of a linking system associated to $\Fc$.

\begin{lemma}
Let $(S,\Fc,\Li)$ be a $p$-local finite group with $\Li$ a centric linking system.
If $\Fc$ is constrained, then $|\Li|$ is a classifying space of a model $G$ of $\Fc$.
\end{lemma}

\begin{proof}
By Proposition \ref{model}, we can assume that $\Li=\Li_S^c(G)$.
Set \[\h=\left\lbrace P\in\Ob(\Li)\; |\;P\geq O_p(G)\right\rbrace\]
and let $\Li^\h$ be the full subcategory of $\Li$ with set of objects $\h$. By \cite[Proposition 1.6]{5a1}, $\h$ contains all $\Fc$-centric and $\Fc$-radical subgroups. Thus, by Proposition \ref{object linking system}, $\Li^\h$ is a linking system associated to $\Fc$ and, by Theorem \ref{differentLisimeq}, $|\Li^\h|\cong|\Li|$.

It remains to prove that that $|\Li^\h|\cong BG$. For that purpose, consider, the following functor.
 \[\begin{array}{rrcl}
F\colon & \Li^\h &\longrightarrow & \Li^{\{O_p(G)\}}\\
  &  P\in\Li^\h  &\longmapsto & O_p(G) \\
  &  g\in T_G(P,Q) &\longmapsto & g\in N_G(O_p(G))=G
  \end{array}\]
It gives us a retraction by deformation of $|\Li^\h|$ on the geometric realization of the full subcategory of $\Li$ with unique object $O_p(G)\leq S$.
 As $\Aut_\Li(O_p(G))=N_G(O_p(G))=G$, this last category is $\Bc(G)$. in particular, its geometric realization is a classifying space of $G$. 
\end{proof}

\begin{proposition}
 Let $G$ be a finite group and $S$ a Sylow $p$-subgroup of $G$.
 If we have $C_G(O_p(G))\leq O_p(G)$, then, for every $\Zp[G]$-module $M$,
 the inclusion of $S$ in $G$ induces a natural isomorphism
 \[H^*(G,M)\cong H^*(\Fc^c_S(G),M).\]
\end{proposition}

\begin{proof}
Let $(S,\Fc,\Li)=(S,\Fc_S(G),\Li_S^c(G))$. By assumption, $\Fc_S(G)$ is constrained and $G$ is a model of $\Fc_S(G)$. From Cartan-Eilenberg Theorem, we know that 
\[\xymatrix{Res_S^G\colon H^*(G,M)\ar[r] & H^*(S,M)}\]
is injective and that $\im(\Res_S^G)=\limproj{\T_S(G)} H^*(-,M)$.
Moreover,  
\[H^*(\Fc^c,M)=\limproj{\Fc^c} H^*(-,M)=\limproj{\Li} H^*(-,M)=\limproj{\T_S^c(G)} H^*(-,M)\geq \limproj{\T_S(G)} H^*(-,M).\]
Thus, it remains to prove that $\limproj{\T_S^c(G)} H^*(-,M)\leq \limproj{\T_S(G)} H^*(-,M)$.

Let then $x\in H^*(\Fc^c,M)=\limproj{\T^c_S(G)} H^*(-,M)$.
For $P\leq S$ and $g\in N_G(P,S)$  we have, in $\T_S(G)$, the following commutative diagram
  \[
  \xymatrix{
  PO_p(G)\ar[r]^-g &  gPg^{-1}O_p(G)\\
  P\ar[u]^{e} \ar[r]^-g &  gPg^{-1} \ar[u]_{e}}
 \]
 where $e$ is the trivial element of $G$. Hence, as the top part of the diagram is in $T^c_S(G)$ and $x\in\limproj{\T_S^c(G)} H^*(-,M)$, 
\begin{align*}
c_g^*\circ \Res_{gPg^{-1}}^S(x) &=\Res_P^{PO^p(G)}\circ c_g^*\circ \Res_{gPg^{-1}O_p(G)}^S(x)\\
&=\Res_P^{PO^p(G)}\circ \Res_{PO_p(G)}^S(x)\\
&=\Res_P^S(x).
\end{align*}
Thus $x\in\limproj{\T_S(G)} H^*(-,M)$ and this complete the proof.
\end{proof}

\begin{corollary}\label{isoconstrained}
 Let $(S,\Fc,\Li)$ be a $p$-local finite group.
 If $\Fc$ is constrained and $M$ is a $\Zp[\pi_\Li]$-module, then $\delta_S$ induces a natural isomorphism,
 \[
  H^*(|\Li|,M)\cong H^*(\Fc^c,M).
 \]
\end{corollary}

\section{Actions factoring through a $p'$-group}\label{sec p'}

In this section, for each $p$-local finite group $(S,\Fc,\Li)$ we will assume that \textbf{$\Li$ is a centric linking system}.

\begin{lemma}
 Let $(S,\Fc,\Li)$ be a $p$-local finite group and $(S,O^{p'}(\Fc),O^{p'}(\Li))$ its minimal $p$-local subgroup of index prime to $p$.
If $M$ is a $\Zp[\pi_\Li]$-module, then the inclusion $O^{p'}(\Li)\subseteq\Li$ induces the following isomorphism,
 \[
H^*(|\Li|,M)\cong H^*(|O^{p'}(\Li)|,M)^{\pi_{\Li}/O^{p'}(\pi_\Li)}.                                         
 \]

\end{lemma}

\begin{proof}
By Theorem \ref{p'cover}, $|O^{p'}(\Li)|$ is, up to homotopy, 
a covering space of $|\Li|$ with fundamental group $O^{p'}(\pi_\Li)\trianglelefteq\pi_\Li$. It gives us a fibration sequence 
\[|O^{p'}(\Li)|\rightarrow |\Li|\rightarrow  B\left(\pi_\Li/O^{p'}(\pi_\Li)\right).\]
Consider then the Serre spectral sequence associated
\[H^{s+t}(|\Li|,M)\Leftarrow H^s\left(\pi_\Li/O^{p'}(\pi_\Li),H^t(|O^{p'}(\Li)|,M)\right).\]
$M$ is a $\Zp$-module, thus $H^q(|O^{p'}(\Li)|,M)$ is also a $\Zp$-module.
As $\pi_\Li/O^{p'}(\pi_\Li)$ is a $p'$-group, the $E_2$-page is concentrated in the first column with terms \[H^t(|O^{p'}(\Li)|,M)^{\pi_\Li/O^{p'}(\pi_\Li)}.\]
Thus the spectral sequence collapses on the $E_2$-page and the Lemma follows. 
\end{proof}

\begin{lemma}
 Let $(S,\Fc,\Li)$ be a $p$-local finite group and $(S,O^{p'}(\Fc),O^{p'}(\Li))$ its minimal $p$-local subgroup of index prime to $p$.
If $M$ is a $\Zp[\pi_\Li]$-module, then
\[
 H^*(\Fc^c,M)= H^*(O^{p'}(\Fc)^c,M)^{\Aut_\Fc(S)/\Aut_{O^{p'}(\Fc)}(S)}.
\]
\end{lemma}

\begin{proof}
Notice first that, by Proposition \ref{p and p' sub prop}, $\Ob(O^{p'}(\Fc)^c)=\Ob(\Fc^c)$. Hence, we are working with the same underlying set of objects. Thus, by definition, $H^*(\Fc^c,M)\subseteq H^*(O^{p'}(\Fc)^c,M)^{\Aut_\Fc(S)/\Aut_{O^{p'}(\Fc)}(S)}$. 
On the other hand, by Proposition \ref{p'generate}, we have $\Fc=\langle O^{p'}(\Fc),\Aut_\Fc(S)\rangle$ which gives the converse inclusion.
\end{proof}

\begin{theorem}\label{main}
Let $(S,\Fc,\Li)$ be a $p$-local finite group and $(S,O^{p'}(\Fc),O^{p'}(\Li))$ its minimal $p$-local subgroup of index prime to $p$.
If $M$ is a $\Zp[\pi_\Li]$-module and if the inclusion $\delta_S$ induces an isomorphism 
\[H^*(|O^{p'}(\Li)|,M)\cong H^*(O^{p'}(\Fc)^c,M),\] 
then $\delta_S$ induces an isomorphism
\[
 H^*(|\Li|,M)\cong H^*(\Fc^c,M).
\]
\end{theorem}

\begin{proof}
Recall that, by \ref{p'cover}, $\pi_1(|O^{p'}(\Li)|)=O^{p'}(\pi_{\Li})$. Then we have the following commutative diagram. 
\[\xymatrix{ \Bc (S) \ar[r]^{\delta_S}\ar[rd]_{\delta_S} & \Li \ar[r]^-{\omega}         & \Bc \pi_\Li \ar[r]          & \Bc(\Aut(M))\\
                                                      & O^{p'}(\Li) \ar[r]^-{\omega}\ar[u] & \Bc\left(O^{p'}(\pi_{\Li})\right)\ar[u]\ar[ur]&}
\]
Moreover, by Proposition \ref{p'generate} and Theorem \ref{p'cover}, the projection $\xymatrix{\pi\colon\Li\ar[r] & \Fc}$ induces an isomorphism \[\pi_\Li/O^{p'}(\pi_{\Li})\cong \pi_1(|\Fc^c|)\cong \Aut_\Fc(S)/\Aut_{O^{p'}(\Fc)}(S).\]
Then, by the two previous Lemmas, we obtain
\begin{align*}
H^*(|\Li|,M) &\cong H^*(|O^{p'}(\Li)|,M)^{\pi_{\Li}/O^{p'}(\pi_{\Li})}\\
              &\cong \left(\lim\limits_{\substack{\longleftarrow \\ O^{p'}(\Fc)^c}}H^*(-,M)\right)^{\Aut_\Fc(S)/\Aut_{O^{p'}(\Fc)}(S)}\\
              &\cong  H^*(\Fc^c,M).
\end{align*}
For the second isomorphism, we have to be careful with respect to the action of $\pi_\Li$ on the left side of the isomorphism and $\Aut_\Fc(S)$ on the right side.
In fact here, by Definition \ref{fcstable} of  $\Fc^c$-stable elements, we can see it on the chain level.
The map $\xymatrix{\delta_S^*\colon H^*(|O^{p'}(\Li)|,M)\ar[r] & H^*(S,M)} $, induced by $\xymatrix{\delta_S\colon BS\ar[r] & |O^{p'}(\Li)|}$, gives on the chain level,
\[
 \xymatrix@R=1mm{\Hom_{\Zp[S]}\left(C_*\left(\widetilde{|O^{p'}(\Li)|}\right),M\right)\ar[r] & \Hom_{\Zp[\pi_{O^{p'}(\Li)}]}(C_*(|\Ec(S)|),M)\\
                        f\ar@{|->}[r] & f|_{C_*(|\Ec(S)|)}}
\]
where $\Ec(S)$ is defined as the category 
with set of object $S$ and for each $(s,s')\in S$, $\Mor_{\Ec(S)}(s,s')=\{\varphi_{s,s'}\}$ (in particular $|\Ec(S)|$ is a universal covering space of $BS$). 
Then, for $\varphi\in\Aut_S(\Fc)$, if we choose a lift $\widetilde{\varphi}\in\Aut_\Li(S)$, $\varphi$ acts on the left side by
\[
\xymatrix{f\ar@{|->}[r] & \omega(\widetilde{\varphi}^{-1})f\omega(\widetilde{\varphi})},
\]
and on the right side by,
\[
 \xymatrix{f\ar@{|->}[r] & \omega(\widetilde{\varphi})^{-1}f\circ \varphi^*}.
\]
Finally, the action of $\varphi$ on $\Ec(S)$ corresponds to the action of $\omega(\widetilde{\varphi})$ on $|\Ec(S)|$ (indeed, a lift of $\omega(\widetilde{\varphi})$ in $\widetilde{|O^{p'}(\Li)|}$  joins every vertex $s\in S$ of $|\Ec(S)|$ to the vertex $\varphi(s)$ and similarly for higher simplices).
 Hence, the two actions coincide.
\end{proof}

\begin{corollary}\label{p'iso}
Let $(S,\Fc,\Li)$ be a $p$-local finite group and $M$ be a $\Zp[\pi_\Li]$-module.
If the action of $\pi_\Li$ on $M$ factors through a $p'$-group then $\delta_S$ induces an isomorphism,
\[
 H^*(|\Li|,M)\cong H^*(\Fc^c,M).
\]
\end{corollary}

\begin{proof}
By Theorem \ref{main}, it is enough to prove that $\delta_S$ induces an isomorphism 
\[H^*(|O^{p'}(\Li)|,M)\cong H^*(O^{p'}(\Fc)^c,M).\]
But, as the action on $M$ factor through a $p'$-group, $\pi_1(|O^{p'}(\Li)|)=O^{p'}(\pi_\Li)$ acts trivially on $M$ and Theorem \ref{trivial} gives the wanted isomorphism. 
\end{proof}

We already know, from a previous article (\cite[Theorem 4.3]{Mo1})
that, if $M$ is a finite $\Zp[\pi_\Li]$-module and the action of $\pi_\Li$ on $M$ factor through a $p$-group, then $\delta_S$ induces an isomorphism
\[
 H^*(|\Li|,M)\cong H^*(\Fc^c,M)
\]
(it is a direct corollary of \cite[Theorem 4.3]{Mo1} because, any action of a $p$-group on an abelian $p$-group is nilpotent). 
Hence, with the same arguments, we get another corollary of Theorem \ref{main}.

\begin{corollary}\label{pnilpotent}
Let $(S,\Fc,\Li)$ be a $p$-local finite group and $M$ be a finite $\Zp[\pi_\Li]$-module.
If the action of $\pi_\Li$ on $M$ factors through an extension of a normal $p$-group by a $p'$-group then $\delta_S$ induces an isomorphism,
\[
 H^*(|\Li|,M)\cong H^*(\Fc^c,M).
\]
\end{corollary}

\section{Realizable fusion systems and actions factoring through a $p$-solvable group}\label{sec psolvable}

Consider here a finite group $G$, $S$ a Sylow $p$-subgroup of $G$ and let $(S,\Fc,\Li)$ be the associated $p$-local finite group with $\Li=\Li_S^c(G)$. Set $\T=\T_S^c(G)$ be the centric transporter category of $G$, $\Li^q=\Li_S^q(G)$ be the quasicentric linking system associated to $G$ and $\T^q=\T_S^q(G)$ be the associated quasicentric transporter category. We also write $\pi_\T=\pi_1(|\T|)$. 

 We have a functor
\[\xymatrix{\rho\colon\T_S(G)\ar[r] &\Bc(G)}\]
which sends each object in the source to the unique object $o_G$ in the target and sends, for every $P,Q\leq S$, $g\in T_G(P,Q)$ to $g\in G=\Mor_\Bc{(G)}(o_G)$. As $|\Bc(G)|=BG$, this induces an homomorphism 
\[\xymatrix{\rho_*\colon\pi_\T\ar[r] & G.}\]

Here for $M$ a $\Zp[\pi_\Li]$-module, with action $\varphi\colon\pi_\Li\rightarrow \Aut(M)$ we will suppose that we have the following commutative diagram for some homomorphism $\overline{\varphi}\colon G\rightarrow\Aut(M)$.
\[\xymatrix{ & \pi_\Li \ar[rd]^\varphi& \\ 
\pi_\T\ar[ru]^{\delta_*} \ar[dr]_{\rho_*} & & \Aut(M)\\
                         & G \ar[ru]_{\overline{\varphi}} & }\]
Then, we can compare the cohomology of $|\Li|$ and the cohomology of $G$ when the action factors through a $p$-solvable group.
The main ingredients that we will use are $p$-local subgroups of index a power of $p$ or prime to $p$. 

The following lemma allows us to compare $H^*(|\Li|,M)$ and $H^*(|\T|,M)$.

\begin{lemma}\label{homofunctor}
Let $G$ be a finite group and $(S,\Fc,\Li)$ be an associated $p$-local finite group. Let $\T=\T_S^{\Ob(\Li)}(G)\subseteq \T^q$ be the transporter category associated to $G$ with set of objects $\Ob(\Li)$.
If $M$ is a $\Zp[\pi_\Li]$-module, then the canonical functor
 $\delta\colon\T\rightarrow\Li$ induces a natural isomorphism $H^*(|\T|,M)\cong H^*(|\Li|,M)$.
\end{lemma}

\begin{proof} 
This is a consequence of \cite[Lemma 1.3]{BLO1} with $\Cc=\T$, $\Cc'=\Li$ and the functor $T\colon\Li^{\text{op}}\rightarrow\Zp\Mod$ which sends each object to $M$,
 and each morphism to its action on $M$. Then $\delta$ induces a natural isomorphism $\limproj{\T}{}^*(M)\cong \limproj{\Li}{}^*(M)$. Then 
 \begin{equation*}
 H^*(|\T|,M)= \limproj{\T}{}^*(M)\cong \limproj{\Li}{}^*(M)=H^*(|\Li|,M)
 \end{equation*}
 Where the first and last equality is just an interpretation in terms of functor cohomology and can be found in \cite[Proposition 3.9]{LR}.
\end{proof}

\begin{theorem}\label{psolvable}
Let $G$ be a finite group, $S$ a Sylow $p$-subgroup of $G$, $\Li=\Li_S^c(G)$ and $\T=\T_S^c(G)$.
Let $M$ be a $\Zp[\pi_\Li]$-module and assume that we have the following commutative diagram.
\[\xymatrix{ & \pi_\Li \ar[rd]^\varphi& \\ 
\pi_\T\ar[ru]^{\delta_*} \ar[dr]_{\rho_*} & & \Aut(M)\\
                         & G \ar[ru]_{\overline{\varphi}} & }\]
If $\rho_*$ is surjective and $\Gamma=\im(\varphi)=\im(\overline{\varphi})$ is $p$-solvable, then $\delta$ and $\rho$ induce natural isomorphisms
 \[H^*(|\Li|,M)\cong H^*(|\T|,M)\cong H^*(G,M).\] 
\end{theorem}

\begin{proof}
By Lemma \ref{homofunctor}, we just have to show that $\rho$ induces a natural isomorphism $H^*(|\T|,M)\cong H^*(G,M)$. 
We prove this by induction on the minimal number $n$ of extensions by $p$-groups or $p'$-groups we need to obtain $\Gamma$.

If $n=0$, $\Gamma=1$ and the action of $\pi_\T$ on $M$ is trivial, then it follows from \cite[Proposition 4.5]{OV1}. Assume that, if $\Gamma$ is obtained by $n$ extensions, the result is true and suppose that $\Gamma$ is obtained with $n+1$ extensions. Consider then the last one
 \[ 0\rightarrow \Gamma_n\rightarrow \Gamma \rightarrow Q\rightarrow 0.\]
Denote $H=\overline{\varphi}_*^{-1}(\Gamma_n)$. 
Thus $(T, \Fc_H, \Li_H)=(S\cap H, \Fc_{S\cap H}(H),\Li_{S\cap H}^c(H))$  is a $p$-local subgroup of $(S,\Fc,\Li)$ of index a power of $p$ or prime to $p$. 

\underline{If $Q$ is a $p'$-group.}
In that case, $(T, \Fc_H, \Li_H)$ is a $p$-local finite subgroup of index prime to $p$ (defined in \ref{p and p' sub}). 
Then $\Ob(\Fc^c)=\Ob(\Fc_H^c)$,  $\T_H=\T_{S\cap H}^c(H)\subset\T$ and, by \cite[Proposition 4.1(d)]{OV1}, 
this inclusion of categories induces, up to homotopy, a covering space of $|\T|$ with covering group $G/H=Q$.
We then have the following commutative diagram with exact rows (here, $\xymatrix{\ar@{->>}[r]&}$ means onto)
\[\xymatrix{  0\ar[rr] && \pi_{\T_H} \ar[rr] \ar@{>>}[dd] \ar[dr] && \pi_\T \ar[rr] \ar@{>>}[dd] \ar[dr] && Q\ar[rr] \ar@{=}[dd] && 0 & \\
               &0\ar[rr] && \Gamma_n \ar[rr]  && \Gamma \ar[rr]  && Q \ar[rr] \ar@{=}[dl]  \ar@{=}[ul]&& 0 \\                             
            0\ar[rr] && H \ar[rr] \ar@{>>}[ur] && G \ar[rr] \ar@{>>}[ur] && Q \ar[rr] &&0 & }\] 
and the following fibration sequences
\[ \xymatrix{ |\T_H|\ar[r] &|\T|\ar[r] & BQ}\]
\[ \xymatrix{ BH \ar[r] & BG\ar[r] & BQ}.\]
Moreover, $\rho$ induces a morphism of fibration sequences between these two.
  
\underline{If $Q$ is a $p$-group.}
In that case, we have to be more careful on the collection of subgroups of $S$ we are working with. As in the case when $Q$ is a $p'$-group we want to apply \cite[Proposition 4.1(d)]{OV1}. This forces us to use the following collection. Let 
\[\h=\left\{P\in\Ob(\Fc^q)\mid P\cap T\in\Ob(\Fc_H^q)\right\}.\]
Since $H\unlhd G$, no element of $T=S\cap H$ is $G$-conjugate to any element of $S\smallsetminus T$. 
Thus, by \cite[Lemma 3.5]{5a2}, for every $P\in\Ob(\Fc^{cr})$, $P\cap T\in\Ob(\Fc_H^c)\subseteq \Ob(\Fc_H^q)$. In particular, $\Ob(\Fc^{cr})\subseteq\h\subseteq\Ob(\Fc^q)$. Hence if $\Li^\h\subseteq \Li^q$ is the full subcategory of $\Li^q$ with set of objects $\h$, by Proposition \ref{object linking system}(b), $\Li^\h$ defines a linking system associated to $\Fc$. On the level of transporter systems, the inclusions $\T\subseteq\T^q\supseteq \T^\h$ induce natural isomorphisms $H^*(|\T^\h|,M)\simeq H^*(|\T^q|,M)\simeq H^*(|\T|,M)$.
Indeed, we have the following commutative diagram.
\[
 \xymatrix{ \T\ar[d]_\delta \ar[r] &\T^q\ar[d]_{\delta} & \T^\h\ar[d]_\delta \ar[l]\\
            \Li \ar[r] &\Li^q & \Li^\h\ar[l]}            
\]
The vertical arrows induce isomorphisms in cohomology by Lemma \ref{homofunctor} and the lower horizontal one induces an isomorphism since, by Theorem \ref{differentLisimeq}, the inclusions of categories $\Li\subseteq \Li^q\supseteq \Li^\h$ induces $|\Li|\simeq|\Li^q|\simeq |\Li^\h|$.
Hence the upper arrows induce isomorphisms $H^*(|\T^\h|,M)\simeq H^*(|\T^q|,M)\simeq H^*(|\T|,M)$.
Finally, By Proposition \ref{p and p' sub prop}, $P\in\Ob(\Fc_H^q)$ if and only if $P\leq T$ and $P\in\h$.
In particular, $\T_H^q\subseteq \T^\h$. Thus we can assume for this part that $\T=\T^\h$ and $\T_H=\T_H^{q}$.
 
We have $\T_H\subseteq \T$ is a transporter system associated to $\Fc_H$ and, by definition of $\h$, the hypotheses of  \cite[Proposition 4.1(d)]{OV1}, are satisfied. Thus this inclusion induces a covering space of $|\T|$ with covering group $G/H=Q$.
Therefore, We have the following diagram with exact rows
\[\xymatrix{  0\ar[rr] && \pi_{\T_H} \ar[rr] \ar@{>>}[dd] \ar[dr] && \pi_\T \ar[rr] \ar@{>>}[dd] \ar[dr] && Q\ar[rr] \ar@{=}[dd] && 0 & \\
               &0\ar[rr] && \Gamma_n \ar[rr]  && \Gamma \ar[rr]  && Q \ar[rr] \ar@{=}[dl]  \ar@{=}[ul]&& 0 \\                             
            0\ar[rr] && H \ar[rr] \ar@{>>}[ur] && G \ar[rr] \ar@{>>}[ur] && Q \ar[rr] &&0 & }\]
and the following fibration sequences
\[ \xymatrix{ |\T_H|\ar[r] &|\T|\ar[r] & BQ}\]
\[ \xymatrix{ BH \ar[r] & BG\ar[r] & BQ}.\]
Moreover, $\rho$ induces a morphism of fibration sequences between these two.

Hence, in both cases, we have the following Serre spectral sequences
\[H^{s+t}(|\T|,M)\Leftarrow H^s(Q,H^t(|\T_H|,M)),\]
\[H^{s+t}(G,M)\Leftarrow H^s(Q,H^t(H,M)),\]
and $\rho$ induces a morphism $\rho^*$ of spectral sequences between these two.
By induction, $\rho^*$ gives an isomorphism on the $E_2$ page and then induces an isomorphism of spectral sequences. 
In particular, $\rho$ induces a natural isomorphism 
\[H^*(|\T|,M)\cong H^*(G,M).\]

The result follows by induction.
\end{proof}

Assume the hypotheses of Theorem \ref{psolvable}.
It remains to compare $H^*(G,M)$ with the $\Fc^c$-stable elements. This is also not obvious and they are not isomorphic in all cases. 
On one hand, by Cartan-Eilenberg Theorem, we have $H^*(G,M)\cong \limproj{\T_S(G)} H^*(-,M)$. On the other hand, we have $H^*(\Fc^c,M)=\limproj{\Li_S^c(G)} H^*(-,M)=\limproj{\T_S^c(G)} H^*(-,M)$. Hence, it remains to compare $ \limproj{\T_S(G)} H^*(-,M)$ and  $\limproj{\T_S^c(G)} H^*(-,M)$. For that we can use a result of Grodal \cite{Gr}.

\begin{definition}\label{M essential}
let $G$ be a finite group, $S\in\Syl_p(G)$ and $M$ be a $\Zp[G]$-module.
Let $K$ be the kernel of $G\rightarrow\Aut(M)$.
A subgroup $P\leq S$ is called \emph{$M$-essential} if
\begin{enumerate}[(i)]
\item the poset of non trivial $p$-subgroup of $N_G(P)/P$ is empty or disconnected,
\item $Z(P)\cap K\in\Syl_p(C_G(P)\cap K)$,
\item $O_p(N_G(P)/(P(C_G(P)\cap K)))=1$.
\end{enumerate}
\end{definition}

The property (ii) looks like the definition of $p$-centric and (iii) looks like the definition of $\Fc$-radical. 
For the property (i), if $P$ is $\Fc$-centric and fully normalized in $\Fc$, it is equivalent to $P=S$ or $P$ is $\Fc$-essential (\cite[Definition I.3.2]{AKO}). 
\begin{theorem}[{\cite[Corollary 10.4]{Gr}}]
Let $G$ be a finite group, $S$ a Sylow $p$-subgroup of $G$ and $M$ a $\Zp[G]$-module.

Let $\h$ be a family of subgroup of $S$ containing $S$ and all the subgroups which are $M$-essential.

Then, the inclusion of $S$ in $G$ induce a natural isomorphism,
\[H^*(G,M)\cong \limproj{\T_S^\h(G)} H^*(-,M).\]
\end{theorem}

From this Theorem and Theorem \ref{psolvable}, we get the following Corollary.

\begin{corollary}\label{psolvable2}
 Let $G$ be a finite group, $S$ a Sylow $p$-subgroup of $G$ and $(S,\Fc,\Li)$ the associated $p$-local finite group.
 Let $M$ be a $\Zp[\pi_\Li]$-module and assume that we have the following commutative diagram,
 \[\xymatrix{ & \pi_\Li \ar[rd]^\varphi & \\
           \pi_\T\ar[ru]^{\delta_*}\ar[rd]_{\rho_*} & & \Aut(M)\\
              & G \ar[ru]_{\overline{\varphi}} &          }\]
that $\rho_*$ is surjective and that $\Gamma:=\im(\varphi)=\im(\overline{\varphi})$. 
If $\Gamma$ is $p$-solvable and and all the $M$-essential subgroups of $S$ are $p$-centric, then $\delta$ and $\rho$ induce natural isomorphisms,
\[H^*(|\Li|,M)\cong H^*(G,M) \cong H^*(\Fc^c,M).\]
\end{corollary}

 We also conjecture that it can be generalized to any abstract $p$-local finite group and any $\Zp[\pi_\Li]$-module with a $p$-solvable action.
 
\begin{conjecture}\label{conjecture}
 Let $(S,\Fc,\Li)$ be a $p$-local finite group and let $M$ be a $\Zp[\pi_1(|\Li|)]$-module.
 If the action of $\pi_1(|\Li|)$ on $M$ is $p$-solvable, then the inclusion of $BS$ in $|\Li|$ induces a natural isomorphism
\[
 \xymatrix{ H^*(|\Li|,M)\ar[r]^\cong & H^*(\Fc^c,M).}
\]  
\end{conjecture}

Corollary \ref{pnilpotent} and Corollary \ref{psolvable2} give good evidence for Conjecture \ref{conjecture} to be true. 

The next section, which is a bit technical, is dedicated to give an example of Conjecture \ref{conjecture} where Corollary \ref{psolvable2} doesn't apply (see Remark \ref{example}).

\section{The $p$-local structure of wreath products by $C_p$: an example for Conjecture \ref{conjecture}}\label{sec example}

Let $G_0$ be a finite group, $S_0$ a Sylow $p$-subgroup of $G_0$ and $(S_0,\Fc_0,\Li_0)$ be the associated $p$-local finite group.  
We are interested in the wreath product $G=G_0\wr C_p$, $S=S_0\wr C_p$ and the associated $p$-local finite group $(S,\Fc,\Li)$.
By \cite[Theorem 5.2 and Remark 5.3]{CL}, we have that $|\Li|\simeq |\Li_0|\wr BC_p:=|\Li_0|^p\times_{C_p} E C_p$  and an extension $(\pi_{\Li_0})^p\rightarrow \pi_\Li\rightarrow C_p$. In addition we have a section $C_p\rightarrow \pi_\Li$ coming from $\ast\wr BC_p\rightarrow |\Li_0|\wr BC_p$ and thus $\pi_\Li=\pi_{\Li_0}\wr C_p$. 

We first give a lemma on strongly $p$-embedded subgroups.
For a finite group $G$, a subgroup $H< G$ is \emph{strongly $p$-embedded}, 
if $p\mid|H|$ and for each $x\in G\setminus H$, $H\cap xHx^{-1}$ has order prime to $p$.

\begin{lemma}\label{substrpemb}
 Let $G$ be a finite group, $G_0\leq G$ a subgroup of index a power of $p$.
 If $G$ contains a strongly $p$-embedded subgroup and $p\mid |G_0|$, then $G_0$ contains a strongly $p$-embedded subgroup.
\end{lemma}

\begin{proof}
Let $H$ be a strongly $p$-embedded subgroup of $G$. 
By \cite[Proposition A.7]{AKO}, $H$ contains a Sylow $p$-subgroup of $G$ so, up to conjugacy, we can choose $H$ such that $H$ contains a Sylow $p$-subgroup of $G_0$.
 Hence $G_0\cap H$ contains a Sylow $p$-subgroup of $G_0$ and $p\mid|G_0\cap H|$. We will show that $G_0\cap H$ is a strongly $p$-embedded subgroup of $G_0$.

As $[G: H]$ is prime to $p$ and $[G:G_0]$ is a power of $p$, $G_0\cap H$ is a proper subgroup of $G_0$.

It remains to show that, for each $x\in G_0\setminus G_0\cap H$, $(G_0\cap H)\cap x(G_0\cap H)x^{-1}$ has order prime to $p$.
But $(G_0\cap H)\cap x(G_0\cap H)x^{-1}\leq H\cap xHx^{-1}$, thus, as $H$ is a strongly $p$-embedded subgroup of $G$, this last
subgroup has order prime to $p$ for every $x\in G \setminus H$. In particular, 
for each $x\in G_0\setminus G_0\cap H$, $(G_0\cap H)\cap x(G_0\cap H)x^{-1}$ has order prime to $p$ and $G_0\cap H$ is a strongly
$p$-embedded subgroup of $G_0$.
\end{proof}

We give also a lemma on $\Fc_1$-essential subgroups for $\Fc_1\subseteq \Fc$ a subsystem of index a power of $p$. 
A proper subgroup $P<S$ is \emph{$\Fc$-essential} if $P$ is $\Fc$-centric and fully normalized in $\Fc$, and if $\Out_\Fc(P)$ contains a strongly $p$-embedded
subgroup.

\begin{lemma}\label{subpessential}
Let $(S,\Fc,\Li)$ be a $p$-local finite group and $(S_1,\Fc_1,\Li_1)$ a $p$-local subgroup of index a power of $p$.
If $P< S_1$ is $\Fc$-essential, then $P$ is $\Fc_1$-conjugate to an $\Fc_1$-essential subgroup and $P$ is $\Fc_1$-essential  if and only if $P$ is fully normalized in $\Fc_1$.
\end{lemma}

\begin{proof}
Let $P< S_1$ be an $\Fc$-essential subgroup. Since $\Fc_1$ is saturated, $P$ is $\Fc_1$-conjugate to a subgroup of $S_1$ fully normalized in $\Fc_1$. If $P$ is $\Fc_1$-essential, it is in particular fully normalized in $\Fc_1$. Thus, it remains to prove that if $P$ is fully normalized in $\Fc_1$, then $P$ is $\Fc_1$-essential. For the remaining, we assume that $P$ is fully normalized in $\Fc_1$ and $\Fc$-essential.

\underline{$P$ is $\Fc_1$-centric:}
As $P$ is $\Fc$-centric, $C_S(Q)=Z(Q)$ for all $Q\in P^\Fc$. In particular, for all $Q\in P^{\Fc_1}\subseteq P^\Fc$, $C_{S_1}(Q)=Z(Q)$ and $P$ is $\Fc_1$-centric.  

\underline{$\Out_{\Fc_1}(P)$ contains a strongly $p$-embedded subgroup:}
Since $P$ is $\Fc$-essential, the group $\Out_\Fc(P)$ contains a strongly $p$-embedded subgroup. As $\Fc_1$ is a subsystem of $\Fc$ of index a power of $p$,
$\Out_{\Fc_1}(P)$ is a subgroup of $\Out_\Fc(P)$ of index a power of $p$. Moreover, as $P$ is a proper subgroup of $S_1$, $P<N_{S_1}(P)$ and, 
as $P$ is $\Fc_1$-centric, every element of $N_{S_1}(P)\smallsetminus Z(P)$ induces a non trivial element in $\Out_{\Fc_1}(P)$.
Hence $p\mid |\Out_{\Fc_1}(P)|$ and, by Lemma \ref{substrpemb}, $\Out_{\Fc_1}(P)$ contains a strongly $p$-embedded subgroup.
 \end{proof}

 We can easily describe the essential subgroups of a product of fusion systems.
 
\begin{lemma}\label{essential product}
 Let $(S_1,\Fc_1,\Li_1)$ and $(S_2,\Fc_2,\Li_2)$ be $p$-local finite groups and set $S=S_1\times S_2$ and $\Fc=\Fc_1\times\Fc_2$.
 The $\Fc$-essential subgroups of $S$ are of the form $Q_1\times S_2$ with $Q_1<S_1$ $\Fc_1$-essential or $S_1\times Q_2$ with $Q_2>S_2$ $\Fc_2$-essential.
\end{lemma}

\begin{proof}
Let $P\leq S$ be a $\Fc$-essential subgroup. By \cite[Proposition I.3.3]{AKO}, $P$ is $\Fc$-centric and $\Fc$-radical.
Thus, by \cite[Lemma 3.1]{AOV}, $P=P_1\times P_2$ with $P_i\leq S_i$ and $P_i$ $\Fc_i$-centric.

Remark also that, if we have two groups $G_1$ and $G_2$ such that $p$ divide $|G_1|$ and $|G_2|$ then $G_1\times G_2$ cannot contain a strongly $p$-embedded subgroup. To see that let $S_i$ be a Sylow $p$-subgroup of $G_i$ and set $H=\langle x\in G\;\mid\; x(S_1\times S_2)x^{-1}\cap S_1\times S_2\neq 1\rangle$. $H$ contains $G_1\times \{0\}$ and $\{0\}\times G_2$ so that $H=G$. Thus, by \cite[Proposition A.7]{AKO}, this implies that $G$ has no strongly $p$-embedded subgroups.

We also have that $\Out_\Fc(P)=\Out_{\Fc_1}(P_1)\times\Out_{\Fc_2}(P_2)$.
Hence, the only possibility for $P$ to be $\Fc$-essential is that $P_1=S_1$ and $P_2$ is $\Fc_2$-essential or the contrary.
\end{proof}

Let $G_0$ be a finite group, $S_0$ a Sylow $p$-subgroup of $G_0$ and $(S_0,\Fc_0,\Li_0)$ be the associated $p$-local finite group. We consider the wreath product $G=G_0\wr C_p$, $S=S_0\wr C_p$ and the associated $p$-local finite group $(S,\Fc,\Li)$.
 Here, for the direct computation, we will take the notation of Alperin and Fong \cite{AF}: an element of $G$ will be represented by permutation matrix corresponding to the powers of $(1,2,\dots, p)$ with entries in $G_0$ and the composition will follow the matrix product with the composition in $G_0$.
 Denote by $c\in G$ the element 
\[ e\otimes P_{(1,2,\dots,p)}=
\begin{pmatrix}      
                  0 & 0      & \cdots & 0      & e     \\
                  e & 0      & \cdots & 0      & 0     \\
             \vdots & \ddots &        & \vdots & \vdots\\
             \vdots &        & \ddots & 0      & 0     \\
                  0 & \cdots & \cdots & e      & 0 
  \end{pmatrix}
\]
where $e$ is the trivial element of $G_0$.
Here, we are interested in the $\Fc$-essential subgroups.
 
 \begin{lemma}\label{essentials}
 Let $P\leq S$ be an $\Fc$-essential subgroup.
 \begin{enumerate}[($E_1$)]
  \item If $P\leq S_0^p$, then either $P=S_0^p$ and $N_G(P)=N_{G_0}(S_0)\wr C_p$
        or $P$ is $\Fc_0^p$-essential and $N_G(P)=N_{G_0^p}(P)$.\label{E0}
  \item If $P\nleq S_0^p$, then $P\cong_\Fc Q\wr C_p$ where $Q$ is $\Fc_0$-essential
  and we have $N_G(P)/P\cong N_{G_0}(Q)/Q$ through the diagonal map $G_0\hookrightarrow G_0^p$.
 \end{enumerate}
\end{lemma}

\begin{proof}
 Let $P\leq S$ be an $\Fc$-essential subgroup.

Assume first that $P\leq S_0^p$. If $P=S_0^p$ a direct calculation gives $N_G(P)=N_{G_0}(S_0)\wr C_p$. 
Else, by Lemma \ref{subpessential}, we know that $P$ is $\Fc_0^p$-conjugate to an $\Fc_0^p$-essential subgroup $Q\leq S_0^p$. By Lemma \ref{essential product} we have $N_G(Q)\leq G_0^p$ and, in particular, $N_G(Q)=N_{G_0^p}(Q)$. Thus, since $P$ is $\Fc_0^p$-conjugate to $Q$, we also have $N_G(P)=N_{G_0^p}(P)$ and, since $P$ is fully normalized in $\Fc$, it is fully normalized in $\Fc_0^p$. Hence, by Lemma \ref{subpessential}, $P$ is $\Fc_0$-essential.
 
Secondly, assume that $P\nleq S_0^p$. 
As all choices of a splitting $C_p\rightarrow G$ are conjugate in $G$, we can assume that $P=\langle P_0,x\rangle$ where $P_0=P\cap S_0^p$ and $x=((x_1,x_2,\dots,x_p),c)$ is such that $x^p\in P_0$. Up to conjugation in $S_0\wr C_p$ we can assume that $x$ is of the form $((a,1,1,\dots,1),c)$ where $a\in N_{S_0}(Q)$ where $Q$ is the projection of $P_0$ on the first factor. 
If we write $P_0^{(i)}$ the projection of $P_0$ on its $i$th factor, as $x$ normalizes $P_0$,
we have that $P_0^{(i)}=P_0^{(j)}$ for all $i,j$ and then $P_0\leq (P_0^{(1)})^p= Q^p$.

Notice also that $N_G(P)=\langle N_{G_0^p}(P),x\rangle$. 
If $g=(g_1,\dots, g_p)\in N_{G_0^p}(P)$, as $g$ normalizes $P\cap G_0^p=P_0$, we have, for all $i$, $g_i\in N_{G_0}(Q)$. 
Moreover, if we denote $h=(h_1,\dots, h_p)=gxg^{-1}x^{-1}\in P_0$, we have, for all $i$, $g_ih_i=g_{i-1}$  (with $g_0=g_p$). Therefore, there is $h'\in Q^p$ such that $g=(g_1,g_1,\dots,g_1).h'\in\left\langle N_{G_0}(Q)\otimes \Id, Q^p\right\rangle\leq N_G(Q^p)$.
Hence, every automorphism $c_g\in \Aut_\Fc(P)$ can be extended to an automorphism of $\langle Q^p,x\rangle$.
As $P$ is $\Fc$ essential, by \cite[Proposition I.3.3]{AKO}, $P=\langle Q^p,x\rangle$.
Now, $x^p\in Q^p$ implies that $a\in Q$ so $P=\langle Q^p,x\rangle=\langle Q^p,c\rangle=Q\wr C_p$.

Finally, direct computations give that 
\[
C_G(P)\cong C_{G_0}(Q)\otimes \Id=\left\{ 
\begin{pmatrix}
 g & 0 & \cdots & 0\\
 0 & g & \ddots & \vdots\\
 \vdots &\ddots &\ddots & 0 \\
 0 & \cdots &  0 & g
\end{pmatrix}
\;;\; g\in C_{G_0}(Q)\right\}
\]
and 
\[N_G(P)/P\cong N_{G_0}(Q)/Q\otimes \Id\cong N_{G_0}(Q)/Q.\] 
In particular, as $P$ is $p$-centric, $Q$ is $G_0$-centric.
Moreover, as $N_G(P)/P=\Out_\Fc(P)$ contains a strongly $p$-embedded subgroup,  $\Out_{\Fc_0}(Q)=N_{G_0}(Q)/Q$ does as well.
Up to conjugacy, we can also assume that $Q$ is fully normalized in $\Fc_0$ and thus $Q$ is $\Fc_0$-essential.
\end{proof}

Let us now look at some cohomological results.
Recall that for a group $G$, a subgroup $H\leq G$, and $M$ a $\F_p[H]$-module, we define the \textit{induced} and \textit{coinduced} $\F_p[G]$-module by, 
\[\Ind_H^G(M)=\F_p[G]\otimes_{\F_p[H]}M \qquad \coInd_H^G(M)=\Hom_{\F_p[H]}(\F_p[G],M).\]
Recall also that, when the index of $H$ in $G$ is finite, these two $\F_p[G]$-modules are isomorphic (by \cite[Lemma 6.3.4]{We}).

\begin{lemma}\label{topowr}
 Let $X$ be a CW complex and denote by $G$ its fundamental group.
 If $X_0$ is a covering space of $X$ with fundamental group $G_0\trianglelefteq G$ of finite index,
 then, for every $\F_p[G_0]$-module $M$, we have a natural isomorphism of $\F_p[G/G_0]$-modules,
 \[
   H^*(X_0,\Ind_{G_0}^G(M))\cong H^*(X_0,M)\otimes_{\F_p} \F_p[G/G_0].
 \]
 Where, on the right side, $G/G_0$ is only acting by translation on $\F_p[G/G_0]$.
\end{lemma}

\begin{proof}
 This can be easily seen on the chain level.
 Let $\tilde{X}$ be the universal covering space of $X$.
 As $\F_p[G/G_0]$-modules, we have the following
 \begin{align*}
  \Hom_{\F_p[G_0]}( C_*(\tilde{X}),\Ind_{G_0}^G(M))  =\bigoplus_{g\in[G/G_0]} \Hom_{\F_p[G_0]}( C_*(\tilde{X}),g.M) 
 \end{align*}
where the action of $G/G_0$ is permuting the terms in the sum. But, each terms in the sum is isomorphic, as (trivial) $\F_p[G/G_0]$-modules, to $\Hom_{\F_p[G_0]}( C_*(\tilde{X}),M)$. Thus 
\begin{align*}
  \Hom_{\F_p[G_0]}( C_*(\tilde{X}),\Ind_{G_0}^G(M)) \cong \Hom_{\F_p[G_0]}( C_*(\tilde{X}),M)\otimes_{\F_p} \F_p[G/G_0]. 
\end{align*}
 This induces the wanted isomorphism in cohomology.
\end{proof}

\begin{proposition}\label{wriso}
Let $G_0$ be a finite group and $(S_0,\Fc_0,\Li_0)$ be the associated $p$-local finite group.
Consider $G=G_0\wr C_p$, $S=S_0\wr C_p$ a Sylow $p$-subgroup of $G$ and $(S,\Fc,\Li)$ the associated $p$-local finite group.
Let $M$ be a $\F_p[\pi_{\Li_0}]$-module. 

If $\delta_{S_0}$ induce natural isomorphisms  
\[H^*(|\Li_0|,M)\cong H^*((\Fc_0)^c,M),\] 
and
\[H^*(|\Li_0|^p,\coInd^{\pi_{\Li}}_{\pi_{\Li_0}^p}(M^{\otimes p}))\cong H^*((\Fc_0^p)^c,\coInd^{\pi_{\Li}}_{\pi_{\Li_0}^p}(M^{\otimes p}))\]
then $\delta_S$ induces a natural isomorphism
  \[
    H^*(|\Li|,\coInd^{\pi_{\Li}}_{\pi_{\Li_0}^p}(M^{\otimes p}))\cong H^*(\Fc^c,\coInd^{\pi_{\Li}}_{\pi_{\Li_0}^p}(M^{\otimes p})).
  \]
\end{proposition}

\begin{proof}
 Write $N=\coInd^{\pi_{\Li}}_{\pi_{\Li_0}^p}(M^{\otimes p})$ and, for $i\in\{1,2\}$,
 denote by $H^*(\Fc^{E_i},N)$ the stable elements of $H^*(S,N)$ under the full subcategory of $\Fc$ 
with objects $S$ and all the subgroups of $S$ of type $(E_i)$ defined in Lemma \ref{essentials}.

By the Mackey Formula, 
\[\Res_{Q\wr C_p}^{\pi_\Li}\Ind_{\pi_{\Li_0}^p}^{\pi_\Li}=\Ind_{Q^p}^{Q\wr C_p}\Res_{Q^p}^{\pi_{\Li_0}^p}.\]
Thus by Shapiro's Lemma (see for example \cite[Proposition 4.1.3]{Ev}) and the Kunneth Formula, for every $P=Q\wr C_p$ of type $(E_2)$, we have a natural isomorphism $H^*(Q\wr C_p,N)\cong H^*(Q^p,M^{\otimes p})\cong H^*(Q,M)^{\otimes p}$ and, by the computation of normalizers in Lemma \ref{essentials},
\[H^*(Q\wr C_p,N)^{\Aut_\Fc(Q\wr C_p)}\cong (H^*(Q,M)^{\Aut_{\Fc_0}(Q)})^{\otimes p}.\]
Hence, applying this to all the subgroups of type $(E_2)$ and, by naturality of the Shapiro isomorphisms, we have that, 
\[H^*(\Fc^{E_2},N)\cong H^*(\Fc^c_0,M)^{\otimes p}.\]
On the other hand, by \cite[Theorem 5.2 and Remark 5.3]{CL}, $|\Li_0|^p$ has the homotopy type of a covering space of $|\Li|$ with covering group $C_p$. 
 Then, if we denote by $X$ the universal covering space of $|\Li|$ (which is also the universal covering space of $|\Li_0|^p$),
we have the following isomorphism on the chain level (because $\Res$ and $\coInd$ are adjoint functors) 
\[\Hom_{\Zp[\pi_{\Li_0}^p]}(C_*(X),M^{\otimes p})\cong \Hom_{\Zp[\pi_{\Li}]}(C_*(X),N)\]
which is analogue to the Shapiro isomorphism (see \cite[Proposition 4.1.3]{Ev}).
By the Kunneth Formula, it gives us the following isomorphism on cohomology
\[H^*(|\Li_0|,M)^{\otimes p}\cong H^*(|\Li|,N)\]
and give the following commutative diagram
 \[\xymatrix{ H^*(S_0,M)^{\otimes p}\ar[r]^{\cong}\ar[d]_{\left(\delta_{S_0}\right)^*} & H^*(S,N)\ar[d]^{\delta_{S}^*}\\
             H^*(|\Li_0|,M)^{\otimes p}\ar[r]_\cong & H^*(|\Li|,N)}                 
 \]
Thus $\delta_S$ induces an isomorphism 
\[H^*(\Fc^{E_2},N)\cong H^*(\Fc_0^c,M)^{\otimes p}\cong H^*(|\Li_0|,M)^{\otimes p}\cong H^*(|\Li|,N).\]

Secondly, by factoring the Shapiro isomorphism (see \cite[Proposition 4.1.3]{Ev}),
the inclusion of $S_0^p$ in $S$ induces an injection $H^*(S,N)\hookrightarrow H^*(S_0^p,N)$.
Hence
\[H^*(\Fc^{E_1},N)\cong H^*((\Fc^p_0)^c,N)^{C_p}\leq H^*(S_0^p,N).\]
By assumption,  $\delta_{S_0^p}$ induces an isomorphism 
\begin{align*}
H^*((\Fc_0^p)^c,N) \cong H^*(|\Li_0|^p,N).
\end{align*}
Moreover, by Lemma \ref{topowr}, this last term is isomorphic to $H^*(|\Li_0|^p,M^{\otimes p})\otimes \F_p[C_p]$ and, in particular, it is a projective $\F_p[C_p]$-module.

Consider now the Serre spectral sequence associated to the fibration sequence \[\xymatrix{ |\Li_0|^p \ar[r] & |\Li| \ar[r] & B C_p}\] 
with coefficients in $N$.
The $E_2$ page is the following,
\[E_2^{s,t}=H^s(C_p,H^t(|\Li_0|^p,N))\]
 and, by projectivity of $H^t(|\Li_0|^p,N)$, the $E_2$ page is concentrated in the 0th column. 
 Hence,  we have that, $H^*(|\Li_0|^p,N)^{C_p}=E_2^{0,*}\cong H^*(|\Li|,N)$.

In conclusion, \[H^*(\Fc^c,N)= H^*(\Fc^{E_1},N)\cap H^*(\Fc^{E_2},N)\cong H^*(|\Li|,N)\]
and the theorem follows.
\end{proof}

This proposition is a bit technical but we will use it in a specific case. Consider $p=5$, the group $G_0=GL_{20}(F_2)$, 
the wreath product $G=G_0\wr C_5$ and $(S_0,\Fc_0,\Li_0)$ and $(S,\Fc,\Li)$ the associated $5$-local finite groups.
By \cite[Theorem 6.3]{Ru}, we know that $(S_0,\Fc_0,\Li_0)$ admits a $5$-local subgroup of index 4 which is exotic $(S_e,\Fc_e,\Li_e)$ and 
that we have a fibration sequence 
\[
\xymatrix{|\Li_e|\ar[r] & |\Li_0| \ar[r] & BC_4.} 
\]
In particular, we have $\pi_\Li/\pi_{\Li_e}^5=C_4\wr C_5$ and we can be interested in comparing $H^*(|\Li|,N)$ and $H^*(\Fc^c,N)$ for
\[N=\F_5[C_4\wr C_5]=\Ind^{\pi_{\Li}}_{\pi_{\Li_0}^5}(M^{\otimes 5})\cong\coInd^{\pi_{\Li}}_{\pi_{\Li_0}^5}(M^{\otimes 5})\]
(the action factors through a finite group) with $M=\F_5[C_4]$.

By Corollary \ref{p'iso}, we have that $\delta_{S_0}$ and $\delta_{S_0^p}$ induce natural isomorphisms  
\[H^*(|\Li_0|,M)\cong H^*((\Fc_0)^c,M)\text{ and}\] 
\[H^*(|\Li_0|^5,\coInd^{\pi_{\Li}}_{\pi_{\Li_0}^5}(M^{\otimes 5}))\cong H^*((\Fc_0^5)^c,\coInd^{\pi_{\Li}}_{\pi_{\Li_0}^5}(M^{\otimes 5})),\]
(for the second isomorphism, notice that $[\Li_0|^5$ has the homotopy type of a linking system associated to $\Fc_0^5$ by \cite[Proposition 2.17]{CL}).
Hence, all the hypothesis of Proposition \ref{wriso} are satisfied and 
\[H^*(|\Li|,N)\cong H^*(\Fc^c,N).\]

\begin{remark}\label{example}
This gives us an example of  isomorphism between the cohomology of $|\Li|$ and the stable elements when the action factors through a $p$-solvable group which cannot be recovered by a previous result. Notice that, even if the fusion system $\Fc$ is realizable, as $\Fc_e$ is exotic, we cannot find a group $G$ with $S\in\Syl_p(G)$ such that $G$ acts on $M$ in the same way as asked in Section \ref{sec psolvable}.  
This example gives us some additional evidence for Conjecture \ref{conjecture}.
\end{remark}

\newcommand{\etalchar}[1]{$^{#1}$}


\newcommand{\etalchar}[1]{$^{#1}$}
\begin{thebibliography}{BCG{\etalchar{+}}07}

\bibitem[5a1]{5a1}
Carles Broto, Nat{\`a}lia Castellana, Jesper Grodal, Ran Levi, and Bob Oliver.
\newblock Subgroup families controlling {$p$}-local finite groups.
\newblock {\em Proc. London Math. Soc. (3)}, 91(2):325--354, 2005.

\bibitem[5a2]{5a2}
Carles Broto, Nat{\`a}lia Castellana, Jesper Grodal, Ran Levi, and Bob Oliver.
\newblock Extensions of {$p$}-local finite groups.
\newblock {\em Trans. Amer. Math. Soc.}, 359(8):3791--3858, 2007.

\bibitem[AF]{AF}
J.~L. Alperin and P.~Fong.
\newblock Weights for symmetric and general linear groups.
\newblock {\em J. Algebra}, 131(1):2--22, 1990.

\bibitem[AKO]{AKO}
Michael Aschbacher, Radha Kessar, and Bob Oliver.
\newblock {\em Fusion systems in algebra and topology}, volume 391 of {\em
  London Mathematical Society Lecture Note Series}.
\newblock Cambridge University Press, Cambridge, 2011.

\bibitem[AOV]{AOV}
Kasper K.~S. Andersen, Bob Oliver, and Joana Ventura.
\newblock Reduced, tame and exotic fusion systems.
\newblock {\em Proc. Lond. Math. Soc. (3)}, 105(1):87--152, 2012.

\bibitem[BLO1]{BLO1}
Carles Broto, Ran Levi, and Bob Oliver.
\newblock Homotopy equivalences of {$p$}-completed classifying spaces of finite
  groups.
\newblock {\em Invent. Math.}, 151(3):611--664, 2003.

\bibitem[BLO2]{BLO2}
Carles Broto, Ran Levi, and Bob Oliver.
\newblock The homotopy theory of fusion systems.
\newblock {\em J. Amer. Math. Soc.}, 16(4):779--856, 2003.

\bibitem[CE]{CE}
Henri Cartan and Samuel Eilenberg.
\newblock {\em Homological algebra}.
\newblock Princeton Landmarks in Mathematics. Princeton University Press,
  Princeton, NJ, 1999.
\newblock With an appendix by David A. Buchsbaum, Reprint of the 1956 original.

\bibitem[CL]{CL}
Nat{\`a}lia Castellana and Assaf Libman.
\newblock Wreath products and representations of {$p$}-local finite groups.
\newblock {\em Adv. Math.}, 221(4):1302--1344, 2009.

\bibitem[Ch]{Ch}
Andrew Chermak.
\newblock Fusion systems and localities.
\newblock {\em Acta Math.}, 211(1):47--139, 2013.

\bibitem[Ev]{Ev}
Leonard Evens.
\newblock {\em The cohomology of groups}.
\newblock Oxford Mathematical Monographs. The Clarendon Press, Oxford
  University Press, New York, 1991.
\newblock Oxford Science Publications.

\bibitem[Gr]{Gr}
Jesper Grodal.
\newblock Higher limits via subgroup complexes.
\newblock {\em Ann. of Math. (2)}, 155(2):405--457, 2002.

\bibitem[LR]{LR}
Ran Levi and K{\'a}ri Ragnarsson.
\newblock {$p$}-local finite group cohomology.
\newblock {\em Homology Homotopy Appl.}, 13(1):223--257, 2011.

\bibitem[Mo1]{Mo1}
R\'emi Molinier.
\newblock Cohomology with twisted coefficients of the classifying space of a
  fusion system.
\newblock {\em Topology and its Applications}, 212:1 -- 18, 2016.

\bibitem[O4]{O4}
Bob Oliver.
\newblock Extensions of linking systems and fusion systems.
\newblock {\em Trans. Amer. Math. Soc.}, 362(10):5483--5500, 2010.

\bibitem[OV1]{OV1}
Bob Oliver and Joana Ventura.
\newblock Extensions of linking systems with {$p$}-group kernel.
\newblock {\em Math. Ann.}, 338(4):983--1043, 2007.

\bibitem[Ru]{Ru}
Albert Ruiz.
\newblock Exotic normal fusion subsystems of general linear groups.
\newblock {\em J. Lond. Math. Soc. (2)}, 76(1):181--196, 2007.

\bibitem[We]{We}
Charles~A. Weibel.
\newblock {\em An introduction to homological algebra}, volume~38 of {\em
  Cambridge Studies in Advanced Mathematics}.
\newblock Cambridge University Press, Cambridge, 1994.

\end{thebibliography}

\end{document}